\documentclass[12pt]{amsart}   	
\usepackage{geometry}              		
  \usepackage{pgf,tikz}
\usepackage{graphicx}	
\usepackage{amssymb}
\usepackage{amsthm}
\usepackage{amsmath}
\usepackage{mathtools}
\usepackage{ytableau}

\usepackage{hyperref}

\usepackage{lipsum}

\usepackage{caption}
\usepackage{subcaption}

\usepackage{setspace}
\usepackage[T1]{fontenc}

\title[A survey of the work of Ian P. Goulden and David M. Jackson]{Combinatorial and Algebraic Enumeration: a survey of the work of Ian P. Goulden and David M. Jackson}

\author[Foley]{Ang\`ele M. Foley}
\address{Department of Physics and Computer Science\\ Wilfrid Laurier University \\ 75 University Avenue West \\ Waterloo, ON, N2L 3C5 \\ Canada}
\email{afoley@wlu.ca}

\author[Morales]{Alejandro H.\ Morales}
\address{Department of Mathematics and Statistics, 
University of Massachusetts, Amherst\\
710 N. Pleasant Street
Amherst, MA 01003 \\ USA} 
\email{amoralesborr@umass.edu} 

\author[Rattan]{Amarpreet Rattan}
\address{Department of Mathematics, 
Simon Fraser University\\ 
8888 University Dr, Burnaby, BC V5A 1S6 \\ Canada}
\email{rattan@sfu.ca} 

\author[Yeats]{Karen Yeats} 
\address{Department of Combinatorics and Optimization,
University of Waterloo\\
200 University Ave W, Waterloo, ON N2L 3G1\\ Canada}
\email{kayeats@uwaterloo.ca}

\date{\today}		

\usepackage{xcolor}
\DeclareMathOperator{\stat}{stat}
\DeclareMathOperator{\hooks}{hooks}

\let\oldmarginpar\marginpar
\renewcommand\marginpar[1]{\-\oldmarginpar[\raggedleft\tiny #1]%
{\raggedright\footnotesize #1}}

\begin{document}

\numberwithin{equation}{section}

\theoremstyle{definition}
\newtheorem{theorem}{Theorem}[section]
\newtheorem{definition}[theorem]{Definition}
\newtheorem{conjecture}[theorem]{Conjecture}
\newtheorem{problem}[theorem]{Problem}
\newtheorem{lemma}[theorem]{Lemma}
\newtheorem{proposition}[theorem]{Proposition}
\newtheorem{corollary}[theorem]{Corollary}
\newtheorem*{theorem*}{Theorem}
\theoremstyle{remark}
\newtheorem{remark}[theorem]{Remark}
\newtheorem{example}[theorem]{Example}
\newcommand{\preetcomment}[1]{{\color{red} #1}}
\newcommand{\ahmcomment}[1]{{\color{blue} #1}}
\newcommand{\amfcomment}[1]{{\color{magenta} #1}}
\newcommand{\bee}{b}

\newcommand{\Sym}{\mathfrak{S}}
\newcommand{\conj}{\mathcal{C}}
\maketitle

\begin{center}
\dedicatory{\em Dedicated to Ian Goulden and David Jackson for their combined 90 years of insight and inspiration that have shaped the fields of algebraic and combinatorial enumeration.}
\end{center}

\begin{abstract} 
In this survey we discuss some of the significant  contributions of Ian Goulden and David Jackson in the areas of classical enumeration, symmetric functions, factorizations of permutations, and algebraic foundations of quantum field theory.  Through their groundbreaking textbook, {\em Combinatorial Enumeration}, and their  numerous research papers, both together and with their many students, they have had an influence in areas of bioinformatics, mathematical chemistry, algorithmic computer science, and theoretical physics.  Here we review and set in context highlights of their 40 years of collaborative work.
\end{abstract}

\section{Introduction} \label{sec:intro}

\vspace{0.25cm}

\begin{figure}
    \centering
    \includegraphics[width=.8\linewidth]{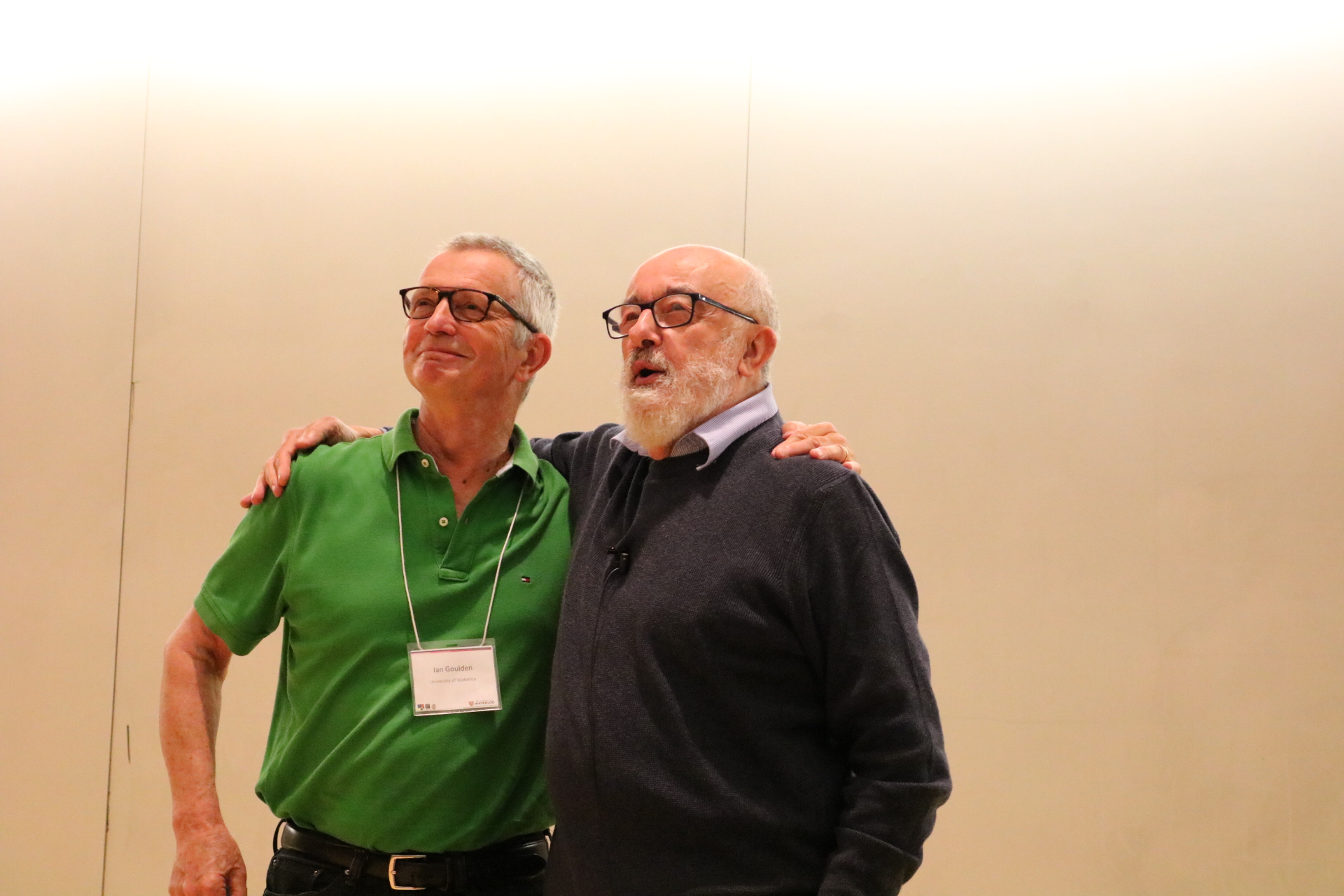}
    \caption{Ian Goulden and David Jackson in May 2022. Photo by Emma Watson.}
    \label{fig:GJ}
\end{figure}

Separately---but most notably together---Ian Goulden and David Jackson (Figure~\ref{fig:GJ}) have made a major impact in the fields of algebraic and combinatorial enumeration. Their contributions run the gamut from writing a groundbreaking textbook and authoritative reference ({\em Combinatorial Enumeration}, Wiley, 1983 \cite{book}), to co-founding with Chris Godsil a well-regarded journal ({\em Journal of Algebraic Combinatorics}, founded in 1992), to publishing over 50 research papers together in their collaboration lasting over 40 years. They have also been invited speakers of the main conference in the field of algebraic combinatorics {\em Formal Power Series and Algebraic Combinatorics (FPSAC)} three times (Goulden in 1995, and Jackson in 1994 and 1996) \cite{fpsac_website}.   Along the way they have supervised numerous graduate students and mentored many more undergraduates and colleagues. This survey article explores some of their key results.

David Jackson graduated from the University of Cambridge in 1970 with a Ph.D. in mathematics. Within three years he had joined the Faculty of Mathematics at the University of Waterloo, where he would spend the rest of his professorial career.  Ian Goulden was one of Jackson's earliest Ph.D. students, graduating from the University of Waterloo in 1979 with a Ph.D. in Statistics.  Shortly thereafter Goulden began his professional career  in the Department of Combinatorics and Optimization at the University of Waterloo.

All of their work spotlights a respect for rigour, a delight in tackling challenging problems, and a drive to understand and explain complex ideas at the interface between algebra and combinatorics. It's a challenging area---the practitioner must be accomplished in two domains of mathematics---but it's a rewarding one. Goulden and Jackson set a high bar. Time and again they have set in motion (or thrown down the gauntlet in) various areas of research. They are responsible for a number of exciting conjectures, some of which we highlight throughout this survey.

Their contributions can be divided into four main threads: classical enumeration, symmetric functions, factorizations of permutations, and quantum field theory, and we explore each of these chronologically and in turn. We have included papers that have garnered significant citations, as determined by {\em Web of Science}, and we have focused primarily on their joint work, or on work that involved their graduate students.

In addition to the work we discuss in this survey, we would like to briefly mention two books.  In 2000, David Jackson and Terry Visentin (Jackson's former Ph.D. student) wrote the book {\em An Atlas of the Smaller Maps in Orientable and Nonorientable Surfaces} \cite{JV_book}. This self-contained volume includes background theory on maps as well as applications, notably to the quadrangulation conjecture (see Section~\ref{sec: maps2factorizations}) and the $b$-conjecture (see Section~\ref{sec: b conjecture}). Comprehensive lists of maps and hypermaps of various types are given (hence the word ``Atlas'' in the title), along with associated tables. The MathSciNet review \cite{mathreview_tutte} is by Bill Tutte, whose own contribution to the theory of maps is well-recognized. 
In 2019, David Jackson and his former postdoc Iain Moffatt wrote the book {\em An Introduction to Quantum and Vassiliev Knot Invariants} \cite{JacksonMoffatt}. This book has three parts: an introduction to combinatorial knot theory, the Reshetikhin-Turaev construction of tangle invariants associated to ribbon Hopf algebras, and an exposition to Vassiliev invariants. The MathSciNet review \cite{mathreview_IainDavid_book}  of this book  indicates that the book does not require a deep knowledge of topology or more advanced topics like tensor categories or $3$-manifolds, making the topic more accessible.

This survey is organized as follows. Section \ref{early-work} discusses the fundamental textbook and early work of Goulden and Jackson. Section \ref{sym-fns} focuses on their work in symmetric functions.  Section \ref{factorizations} explores some of their most important work, that on  factorizations of permutations, including work on the Harer-Zagier formula, maps on surfaces,  Hurwitz numbers, the KP hierarchy, and the $b$-conjecture and matchings-Jack conjecture. Section \ref{quantum} collects some of Jackson's work towards an algebraic foundation for quantum field theory.  

\section{Early work}
\label{early-work}

``Let $R$ be a ring with unity.'' Thus begins the classic text {\em Combinatorial Enumeration} \cite{book},  published originally in 1983.  From the very first sentence we can see an intimation of the goal of injecting rigour into the discipline. This continues for over 500 pages,  
where classical results from authors like MacMahon, Newcomb, and Cayley appear alongside the state-of-the-art results, including, notably, results from their own work on sequences, permutations and cluster decompositions; see, for example, \cite{aleliunis,formalcalculus1,formalcalculus2, formalcalculus3,algebraicmethods,inver-thm}. Over 350 exercises are included, with complete solutions (which always delights the student and researchers alike!).

The strength of the book is the seamless melding of the combinatorial objects with the algebraic methods from formal power series that are used to describe and manipulate them. Thus we have a full discussion of the Lagrange Inversion theorem, with applications and examples, along with both ordinary and exponential generating functions as applied to permutations, sequences, plane trees, rooted planar maps, Ferrers diagrams, labeled trees, and more unusual objects such as Latin rectangles and  $0,1$ matrices. They provide a detailed program on the combinatorics of sequences, including decomposition theorems based on various patterns within the sequence.  The book includes a study on the combinatorics of paths, both those specified by type of step (e.g.\ up-diagonal, down-diagonal, horizontal) and those specified by height (typically embedded in a lattice). This final chapter culminates with a $q$-analogue of the Lagrange Inversion theorem, thus coming full circle.

Goulden and Jackson state in the preface: ``The book is written not only for the combinatorial theorist but also for the mathematician, the physicist, and the computer scientist, in whose fields problems of this type occur.'' It has fulfilled its promise.
A glance through the over 1600 citations on Google Scholar reveals the book has been referenced in the context of software watermarking \cite{watermarking}, ruin theory in economics \cite{lin-willmot},
Bose-Einstein condensation \cite{condensation}, as well, of course, in many areas of mathematics.  
In 1994,  while leafing through preprints at the University of Canterbury, New Zealand, one of the present authors (AMF) was pleased and astonished to discover {\em Combinatorial Enumeration} cited in the context of mathematical biology \cite{steel}.   Their text remains a classic reference and was re-issued as a Dover reprint in 2004 \cite{bookDover}.

Much of their early generating function and sequence work has been subsumed by the book, but to round out this section we single out some of their papers that  laid the groundwork for their later program on factorizations. Thus we have Jackson's papers \cite{jack:1} and  \cite{jacksomeproblems}, which will be discussed at length in Section \ref{factorizations}, elaborating on the idea of constraining the number of cycles in a permutation and counting the result.

We would also like to highlight in particular two of their papers that employ combinatorial bijections.
One of their most cited papers \cite{goul:7} uses a bijection between $m$-tuples of permutations and cacti, each with appropriate constraints.  A cactus is a generalization of a tree and is a connected graph made up of $m$-gons where every edge lies on exactly one of the $m$-gons.  In this paper Goulden and Jackson start by establishing a bijection between pairs of permutations in $\mathfrak{S}_n$ whose product is the canonical long cycle and subject to an external constraint on cycle numbers, and two-coloured plane edge-rooted trees on $n$ edges. This bijection  preserves  cycle distribution in the permutations and colour distribution of vertices in the tree. The more general result preserves similar statistics in the $m$-tuples of permutations and cacti, and is a beautiful result.
A later single author paper  by Goulden \cite{diffoperator} also uses a bijection involving {\em painted permutations} (in which $n-1$ elements---all but the last---are coloured either red or blue) to provide a bijective proof of a specialization of result of Jackson \cite{jacksomeproblems}, namely that the number of ways of writing an $n$-cycle as a product of $m$ transpositions is
$\frac{1}{n!} \sum_{k=0}^{n-1}( \binom{n}{2} - nk )^m \binom{n-1}{k} (-1)^k.$ 
These ideas will see their full fruition in their work on transitive factorizations and Hurwitz numbers (see  Section~\ref{factorizations}).

\section{Symmetric Functions}
\label{sym-fns}

In the mid 1990's and beyond, Goulden and Jackson, both together and individually, have had an impact on symmetric function theory.  Symmetric functions are key objects in algebraic combinatorics, and are important for their connections to representation theory as well as algebraic geometry.  Symmetric functions form an algebra, and the various bases for symmetric functions---Schur functions, power sum symmetric functions, elementary symmetric functions, homogeneous symmetric functions, and more---are  classical objects of study, and some of the most fundamental questions about them concern the positivity of writing the functions of one basis in terms of another. These twin themes of positivity and coefficients appear repeatedly in the work we consider in this section.

One of their most influential endeavours was a pair of papers on immanants.  An immanant of a matrix is a sum over product of matrix entries multiplied by an irreducible character of the symmetric group.  In symbols,
\[
\sum_{\sigma \in \Sym_n} \chi_\lambda(\sigma) a_{1, \sigma(1)} a_{2, \sigma(2)} \ldots a_{n, \sigma(n)},
\]
where $\lambda$ is a partition of $n$, $\chi_\lambda(\sigma)$ is the irreducible character of the symmetric group, indexed by $\lambda$ and evaluated at $\sigma$, and $A=(a_{ij})$ is a matrix.
The determinant (where the character is the alternating sign character) and the permanent  (where the character is $1$ for all permutations) are both special cases of immanants.
Goulden and Jackson derived a relationship between immanants and Schur functions \cite{immanants1}, showing how an immanant could be written as the Schur function with particular arguments.  They further made connections with this interpretation and early work of Littlewood and McMahon.
At the same time Goulden and Jackson studied immanants of Jacobi--Trudi matrices \cite{immanants2}, showing that if the associated partition was a border strip then the immanant of the Jacobi--Trudi matrix had nonnegative coefficients. They further conjectured that the nonnegativity held for all skew partitions.  This conjecture was proved by Curtis Greene \cite{GouldenGreene} and further generalized by Haiman~\cite{Haiman} using {\em Kazhdan--Lusztig theory}, but the influence of the conjecture went beyond its proof. We can trace a direct path from this conjecture to a paper of Stembridge and Stanley \cite{StanStem}, 
and on to the paper of Stanley \cite{Stan} that defines the chromatic symmetric functions and introduces what has become a celebrated $e$-positivity conjecture for a class of chromatic symmetric functions. As explained in \cite{Stan} this $e$-positivity conjecture had its genesis in this domain.

In addition to the immanants work, Goulden and Jackson have another pair of significant papers in symmetric function theory concerning connection coefficients (see Section \ref{factorizations} of this survey for a definition of connection coefficients).
In \cite{GJ-Jack}  they study a relationship between connection coefficients and {\em Jack symmetric functions} that culminates in the {\em Matchings-Jack Conjecture}, which  is still open, although there are many partial results. Full details, including reference for some of these partial results, is found in Section \ref{factorizations}.

The paper \cite{goul:5} contains a direct proof of a result of Macdonald that the connection coefficients for a particular class of symmetric functions are the same as the connection coefficients of the class algebra of the symmetric group, subject to a ``top'' constraint. The key tool they use in their proof is the combinatorial construction from their earlier paper \cite{goul:7}, connecting cacti and factorizations, as well as Lagrange Inversion, making this a perfect example of how algebraic combinatorics can be leveraged to prove results in algebra.  The methods and results from this paper feature in their later work on transitive factorizations.

Individually Goulden has had an impact in symmetric function theory through papers with various co-authors.  
Several of these papers concern the Schur function, a symmetric function defined in terms of partitions.  A partition $\lambda = (\lambda_1, \lambda_2, \ldots, \lambda_k)$ is a non-increasing set of nonnegative integers, e.g.\ $\lambda=(4, 4, 2, 1, 1)$. We can represent a partition pictorially by a Ferrers diagram $F^{\lambda}$ which is a left justified set of boxes with $\lambda_i$ boxes in the $i$th row 
indexed from top to bottom (using the `English' convention), for $1\leq i \leq n$, e.g.\

\begin{equation}\label{x1}
\ytableausetup{smalltableaux}
\ydiagram{4,4,2,1,1}.
\end{equation}

A semistandard tableau $T$ of shape $\lambda$ is a filling of a Ferrers diagram of shape $\lambda$ with positive integers such that the entries weakly increase in rows and strictly increase in columns, e.g.\

\begin{equation}\label{x2}
\ytableausetup{smalltableaux}
\ytableaushort{1113,3455,45,5,6}.
\end{equation}
$T^{\lambda}$ is the set of semistandard tableaux of shape $\lambda$.  The Schur function is a sum of weighted terms coming from weighted semistandard tableaux, i.e.
\begin{equation} \label{Schurfn}
s_\lambda(x) = \sum_{T\in T^{\lambda}} \prod_{(i,j)\in F^{\lambda}} wgt(t_{ij}),
\end{equation}
where $t_{ij}$ is the entry in box $(i,j)$ in row $i$ and  column $j$ of $T$, and where $wgt(t_{ij})$ is $x_k$ if $ t_{ij}=k$.

Given a second partition $\mu$, we can define a skew partition $\lambda/\mu = (\lambda_1 -\mu_1, \lambda_2 - \mu_2, \ldots, \lambda_n - \mu_n)$ with analogous definitions for Ferrers diagrams, semistandard tableaux, and Schur functions, e.g. the Ferrers diagram for $\lambda/\mu= (4, 4, 2, 1, 1 )/(3,1) $ looks like this:
\begin{equation}\label{x3}
\ydiagram{3+1,1+3,2,1,1},
\end{equation}
and the skew Schur function is defined as
\begin{equation} \label{skewSchurfn}
s_{\lambda/\mu}(x) = \sum_{T\in T^{\lambda/\mu}} \prod_{(i,j)\in F^{\lambda/\mu}} wgt(t_{ij}).
\end{equation}

The paper \cite{GouldenGreene} with Curtis Greene establishes the general definition of factorial symmetric functions which extends the Schur function definition in (\ref{skewSchurfn}):

\[s_{\lambda/\mu} (x|a) = \sum_{T\in {T}^{\lambda/\mu}} \prod_{(i,j)\in F^{\lambda/\mu}}(x_{t_{ij}}+a_{t_{ij}+j-i})
\]
using an arbitrary parameter $a$.  Their discovery was simultaneous with that of Macdonald \cite{Mac92}, but provides a different perspective on these functions. The paper \cite{HamelGoulden} with former Ph.D. student Ang\`ele Hamel (Foley) defines the Hamel--Goulden identity, a general framework for constructing determinantal identities of symmetric functions. The Jacobi-Trudi, dual Jacobi-Trudi, and Giambelli identities are all special cases of the Hamel--Goulden identity, and the identity has found use in different contexts, e.g. \cite{RSvW}, \cite{MoralesPakPanova}.  The  paper \cite{GouldenRattan} on Kerov's character polynomials with former Ph.D. student Amarpreet Rattan is a return to polynomials inspired by character theory and it derives an explicit formula for various components of Kerov's character polynomials.  In this context they also introduce what have become known as {\em Goulden--Rattan polynomials} and they further conjecture positivity results for these polynomials. To date the conjecture has not been proved (or disproved) but partial results have been obtained, e.g.\ \cite{marciniak}.

\section{Factorizations of permutations}
\label{factorizations}

This section concerns the extensive work by Goulden and Jackson in regards to
permutation factorization problems, with a focus on the transitive
factorization problem.  Here we let $\Sym_n$ be the symmetric group on $n$
symbols, which can usually be assumed to be $[n]:=\{1, \ldots, n\}$.  We write $\beta \vdash n$ to indicate $\beta$ is a partition of $n$. If $\alpha$ has $m_i$ parts of size $i$, set $\mathrm{aut}(\alpha) = \prod_i m_i!$.  The length or number of parts of $\beta$ is denoted by $l(\beta)$. By abuse of notation, given a permutation $\pi$, we denote by $l(\pi)$ the number of cycles of $\pi$. We denote by $\conj_\beta$ the conjugacy class of $\Sym_n$ of permutations with cycle
type $\beta$.  We call the special class $\beta = (n)$, the class of 
$n$-cycles in $\Sym_n$, \emph{full cycles}.   Permutation factorization 
problems involve finding the number of
ways of writing a permutation (the \emph{target}) as a product of others (the
\emph{factors}), subject to constraints.  We discuss three different kinds of 
constraints here.  The first two are local constraints on each factor, while 
the third is a global constraint on all the factors.

\subsection{Factorizations of the long cycle} \label{sec:long cycle} Two seminal papers of Jackson involve constraining the number of cycles in a 
factor. 
For example, Jackson \cite{jack:1} shows that in $\Sym_{pn}$, for a full cycle 
$\pi$ and a positive number $k$, the number $e_k^{(p)}(n)$ of ordered pairs  
$(\sigma, \tau)$ such that $\pi = \sigma \tau$, with $\sigma$ having $n$ cycles of size $p$ and $\tau$ having $k$ cycles, is given in terms of {\em Stirling numbers} of the first and second kind, $s(a,b)$ and $S(a,b)$, respectively:
\begin{equation} \label{eq:factorizations and Stirling}
e_k^{(p)}(n) = \frac{1}{(1+pn)p^{n+k}} \sum_{m=n+k}^{pn+1} p^m \binom{m}{k} s(pn+1,m)S(m-k,n).
\end{equation}
In particular $e^{(p+1)}_{pn+1}(n)=\frac{1}{pn+1}\binom{(p+1)n}{n}$, which is a {\em Fuss--Catalan number} \cite[A14]{StanleyCat}.
This result was one of the first of many by 
Jackson and others in this area.  Another example is in 
\cite{jacksomeproblems}, where Jackson proves the following remarkable formula: 
 if $\gamma(k,m;n)$ is the number of ways of writing a fixed $n$-cycle
as the product of two permutations with $k$ and $m$ cycles, respectively, then
\begin{equation}\label{eq:jackbinom}
 \sum_{k,m \geq 1} \gamma(k,m;n) x^k y^m = n! \sum_{p,q \geq 1} \binom{n-1}{p-1,q-1} \binom{x}{p} \binom{y}{q},
\end{equation}
where $\binom{x}{k}:=x(x-1)\cdots (x-k+1)/k!$.
Jackson in fact proves a more general formula for factorizations of a
$n$-cycle into any fixed number of factors, where each factor has a
specified number of cycles: if $\gamma(m_1,m_2,\ldots,m_k;n)$ is the
number of ways of writing a fixed $n$-cycle as the product of $k$
permutations with $m_1,\ldots,m_k$ cycles, respectively, then
\begin{equation} \label{eq:jackpfactors}
\sum_{k,m \geq 1} \gamma(m_1,\ldots,m_p;n) x_1^{m_1}\cdots x_k^{m_k}
\,=\, (n!)^{k-1} \sum_{p_1,\ldots,p_k} M^n_{p_1,\ldots,p_k}  \binom{x_1}{p_1}\cdots \binom{x_k}{p_k},  
\end{equation}
where $M^n_{p_1,\ldots,p_k}$ is the number of tuples
$(S_1,\ldots,S_n)$ of proper subsets of $[k]$ such that $p_j$
of the sets contain $j$. It is open to find a simple combinatorial proof of this result. The tool used to prove these results is using fundamental results in group theory, going back to work of Frobenius \cite{Frobenius}, to write this number of factorizations in terms of  {\em irreducible characters} of the symmetric group. See \cite{SV} for a bijective proof of the case $k=2$, and \cite{BM1,BM2} for a tour de force combinatorial proof for all $k$. 

Jackson also gives the following exponential and
ordinary generating functions for the number $t(n,k)$ of
factorizations of a fixed $n$-cycle into $k$ transpositions: 
  \begin{align}
\sum_{k \geq n-1} t(n,k) \frac{x^{k}}{k!} &=
                                                 \frac{e^{x\binom{n}{2}}}{n!}(1-e^{-xn})^{n-1},
                                                 \label{gf exp}\\
    \sum_{\ell \geq n-1} t(n,\ell) x^{\ell} &= n^{n-2} x^{n-1}
                                              \prod_{k=0}^{n-1}
                                              \left(1 -
                                              xn\left(\frac{n-1}{2}-k\right)\right)^{-1}
                                              \label{gf ord}.     
  \end{align}

  The leading term of the generating functions is 
  \begin{equation}
      \label{eq:Cayley}
  t(n,n-1)=n^{n-2},
  \end{equation}
  a result known to Hurwitz \cite{hurwitz}. By Cayley's formula, this number also
counts the number of
trees with vertices $\{1,\ldots,n\}$. A double counting proof of this result relating factorizations
and trees was given by D\'enes \cite{denes}.  Bijective proofs were given by
Moszkowski \cite{mosz}, Goulden and Pepper \cite{gouldenpepper} and Goulden
and Yong \cite{goul:3}. 

These minimal transposition factorizations of the long cycle also have an important poset interpretation by Biane \cite{Biane02} as {\em maximal chains} in the {\em lattice of non-crossing partitions}.

\begin{remark}
The bijection of Goulden and Yong was used by Feray and Kortchemski \cite{FK2} to study the {\em trajectories} of the elements in $[n]$ viewed as finitely many points on a line in a random typical minimal transposition factorization of the long cycle. The same authors in \cite{FK1} studied the shape of a random typical minimal transposition factorizations of the long cycle by viewing each transposition as a chord in the unit disk. 
\end{remark}

These generating functions spurred generalizations to other groups. 
Since the symmetric group is a {\em reflection group} and
transpositions are {\em reflections}, Equations \eqref{gf
  exp} and \eqref{gf ord} have been generalized to other groups. Chapuy--Stump \cite{CC}
gave a uniform generalization of \eqref{gf exp}  to {\em complex reflection
  groups} $W$ by studying the number $t(W,k)$ of factorizations of a {\em Coxeter element} (an
analogue of the $n$-cycle) into $k$ reflections. Their resulting
generating function is
\begin{equation} \label{eq:CSgf}
\sum_{k \geq 0} t(W,k)\frac{x^k}{k!} \,=\,  \frac{e^{xR}}{|W|}(1-e^{-xh})^{n},
\end{equation}
where $n$ is the rank of the group $W$, $R$ is the number of reflections,
and $h$ is the order of a Coxeter element. The leading term of
this generating function is $n!\cdot  h^n/|W|$, a known generalization of Cayley's formula \eqref{eq:Cayley} for complex reflection groups (indeed,  the symmetric group $W=\Sym_{N}$ has rank $n=N-1$ and Coxeter number $h=N$ and so the formula gives $N^{N-2}=t(N,N-1)$). This generalized formula 
is due to Deligne for the case of real reflection
groups and to Bessis \cite{Bessis_Annals} for the other cases. The original proof of \eqref{eq:CSgf} by
Chapuy--Stump was case-by-case using a character theory approach. Soon
after, Michel \cite{JM} gave a uniform proof for Weyl groups by
using the geometry of Deligne--Lusztig varieties to combine characters
and Douvropoulos \cite{D1} gave a different uniform proof for the
complex reflection groups  using Galois
representations to combine characters. 

The ordinary generating function \eqref{gf ord} has been generalized
to the group $GL_n(\mathbb{F}_q)$ of $n\times n$ invertible matrices
over the finite field $\mathbb{F}_q$ by Lewis--Reiner--Stanton \cite{LRS}, also
using the character theory approach. In
this case, the analogue of the $n$-cycle is a {\em Singer cycle}, the
image of a generator of the $(q^n-1)$-cyclic group
$\mathbb{F}_{q^n}^{\times}$ embedded in $GL_n(\mathbb{F}_q)$. The
leading term in this case is $(q^n-1)^{n-1}$ counting the number of
factorizations of a fixed Singer cycle into $n$ reflections (matrices
that fix a hyperplane). No combinatorial proof of this result is currently
known.

Lastly, \eqref{eq:jackpfactors} has analogues in  complex reflection groups in \cite{LM1} and for $GL_n(\mathbb{F}_q)$ in \cite{LM2}.

\subsection{Harer--Zagier formula}

A related formula to \eqref{eq:jackbinom}, the celebrated formula of Harer and Zagier \cite{harerzagier}, is given by
\begin{equation}\label{eq:harer}
    \sum_{k \geq 1} a_{p,k} x^k = (2p-1)!! \sum_{k \geq 1} \binom{p}{k-1} \binom{x}{k},
\end{equation}
where $a_{p,k}$ is the number of ways to write a fixed full cycle in $\Sym_{2p}$ 
as the product of a fixed point free involution and a permutation with $k$ 
cycles\footnote{Note that $a_{p,k}=e^{(2)}_k(p)$ from \eqref{eq:factorizations and Stirling}.}.  The original context for the formula of Harer and Zagier was in 
finding the Euler characteristic of certain moduli spaces.  Naturally, their 
formula begged for a straightforward combinatorial proof.  Several 
mathematicians \cite{IZHZ,KeHZ,KHZ,PHZ}, including Jackson \cite{JHZ}, offered alternative proofs of \eqref{eq:harer}, but the first 
direct combinatorial proof was by Goulden and Nica \cite{goulnica:1},
related to an ingenious connection to the {\em BEST theorem} (see \cite[Thm. 5.6.2]{EC2}) by Lass \cite{Lass}. This theorem was incorporated into a bijection by Bernardi \cite{OB} between {\em unicellular maps} with coloured vertices with exactly $k$ colours and maps with $k$ vertices on an {\em orientable surface} and with a marked spanning tree. Soon after, Chapuy \cite{GC} and Chapuy--F\'eray--Fusy \cite{CFF} gave bijections proving recurrences for the numbers $a_{p,k}$ related to a recurrence of Harer--Zagier for these numbers.

\subsection{Factorizations and cacti formula}

A second type of constraint is to specify the conjugacy class from which factors may 
come.  More precisely, 
the general problem is for a permutation $\rho$ and partitions $\beta_1, 
\ldots, \beta_k 
\vdash n$, to find the number of ordered tuples $(\tau_1, \ldots 
\tau_k)$ of permutations in $\Sym_n$ such that 
\begin{equation}\label{eq:notrans}
\begin{aligned}
	a)\;\;\;\;\;\;&  \rho = \tau_1 \cdots \tau_k\\
	b)\;\;\;\;\;\;&  \tau_i \in \conj_{\beta_i},  \text{ for } i=1,\ldots,k
\end{aligned}
\hspace{9cm}
\end{equation}
By symmetry, this number  only depends on the 
conjugacy class $\alpha$ of $\rho$, so we denote the number of tuples by 
$c^{\alpha}_{\beta_1, \ldots, \beta_k}$.  The numbers $c^{\alpha}_{\beta_1, 
\ldots, \beta_k}$ can be viewed as follows.
Let $K_\gamma = \sum_{\sigma \in \conj_{\gamma}} \sigma$.  Then $\{K_\gamma : 
\gamma \vdash n\}$ is basis for the centre of the group algebra of the 
symmetric group on $n$ symbols, and
\begin{equation}\label{eq:concoeff}
 c^{\alpha}_{\beta_1, \ldots, \beta_k} = [K_\alpha] K_{\beta_1} \cdots 
K_{\beta_k},
\end{equation}
where $[A] B$ denotes the coefficient of $A$ in $B$.  
The numbers $c^{\alpha}_{\beta_1, \ldots, \beta_k}$ are called 
\emph{connection coefficients} of $\Sym_n$. With this algebraic 
interpretation of connection coefficients, one can use the same character theory approach mentioned in Section~\ref{sec:long cycle} to express them in terms of the irreducible characters of $\Sym_n$.  Often, however, these 
expressions are intractable for the purposes of finding precise enumerative 
results, except in limited circumstances.

Early work on this topic is due to Walkup \cite{walkup}, who 
considers the case $k=2$ in \eqref{eq:concoeff} and $\beta_1, \beta_2 = (n)$ 
for arbitrary $\alpha$;  that is, he determines $c^{\alpha}_{(n), (n)}$. 
Stanley \cite{stanncyc} used characters to solve the same 
problem but for arbitrary $k$ and each $\beta_i = (n)$ and also proved a 
conjecture found in \cite{walkup}.  As pointed out later by Stanley \cite{stanley_publications}, Goulden and Jackson would develop this character approach "much more extensively". 

Goulden and Jackson made significant discoveries here.  For  partitions 
$\gamma, \alpha, \beta \vdash n$, elementary results on permutation products 
give that if the connection coefficient $c^{\gamma}_{\alpha, \beta}$ is 
non-zero, then $l(\alpha) + l(\beta) \leq n + l(\gamma)$.  When 
$\gamma, \alpha$ and $\beta$ satisfy this equation with equality, we call the 
connection coefficients \emph{top}.

B\'edard and Goupil first determined $c^{(n)}_{\alpha, \beta}$ for top 
connection coefficients (so $l(\alpha) + l(\beta) = n+1$) with an 
inductive argument. Goulden and Jackson \cite{goul:7} substantially generalized
this result to the case $c^{(n)}_{\beta_1, \ldots, \beta_k}$, where the
$\beta_i$ satisfy the generalized top condition $\sum_{i=1}^k l(\beta_i) = (k-1)n + 1$:
\begin{equation} \label{eq: cacti}
c^{(n)}_{\beta_1, \ldots, \beta_k} \,=\, n^{k-1} \prod_{i=1}^k \frac{(l(\beta_i)-1)!}{\mathrm{aut}(\beta_i)}.
\end{equation}
Their construction in general involves cacti and Lagrange Inversion. The $k=2$ case of their construction is
especially simple and elegant as it involves plane trees. Later, Goupil and Schaeffer \cite{GS} and Poulalhon and Schaeffer \cite{PS} found positive sum formulas, as opposed to the alternating formulas using characters, but with (exponentially) many terms for all coefficients $c^{(n)}_{\alpha,\beta}$ and $c^{(n)}_{\beta_1,\ldots,\beta_k}$, respectively.

Goulden and Jackson \cite{goul:5} also show that the top coefficients appear as
the structure constants of certain symmetric functions (they credit unpublished
notes of Macdonald with this discovery, but supply the first direct proof).
They also give a number of explicit enumerative results;  for example, they give
the number of minimal factorizations of a full $n$-cycle into $(k+1)$-cycles; that
is, in \eqref{eq:notrans}, we have $\alpha = (n)$ and $\beta_i = (k+1,
1^{n-k-1})$.  This is a generalization of the case $k=1$ originally by D\'enes
\cite{denes}.

Also in \cite{goul:5}, they find the number of
\emph{inequivalent} factorizations of a full cycle into $(k+1)$-cycles:  two
factorizations are \emph{equivalent} if one can be transformed into the other
by a sequence of swapping commuting cycles.  This was done
through the use of the \emph{commutation monoid} of Cartier and Foata \cite{cartfoat}.

\subsection{Factorizations and Maps on surfaces} \label{sec: maps2factorizations}

A {\em rooted map} is a graph embedded on a (locally) orientable surface so that all the faces are homeomorphic to a disc, with one vertex distinguished, called the {\em root vertex}, and one edge incident to the root vertex distinguished, called the {\em root edge}. Important families of maps include {\em planar maps} (e.g. see \cite{Scha_survey}), $3$-face and $4$-face regular maps  called {\em triangulations} and {\em quadrangulations}, respectively, {\em bipartite maps}, and {\em constellations} (e.g. see \cite{LZbook}). We sketch two ways that permutation factorizations can be encoded with maps. These constructions can be traced to Tutte \cite{Tutte84}.

\begin{figure}
\begin{subfigure}[b]{.4\linewidth}
\centering \includegraphics{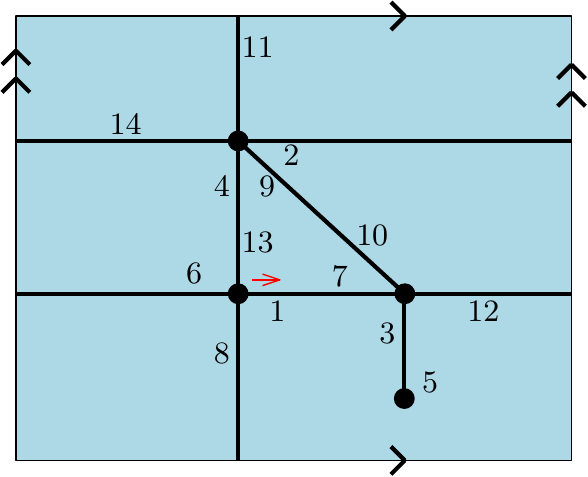}
\caption{}
\label{fig:map}
\end{subfigure}\qquad
\begin{subfigure}[b]{.4\linewidth}
\centering \includegraphics{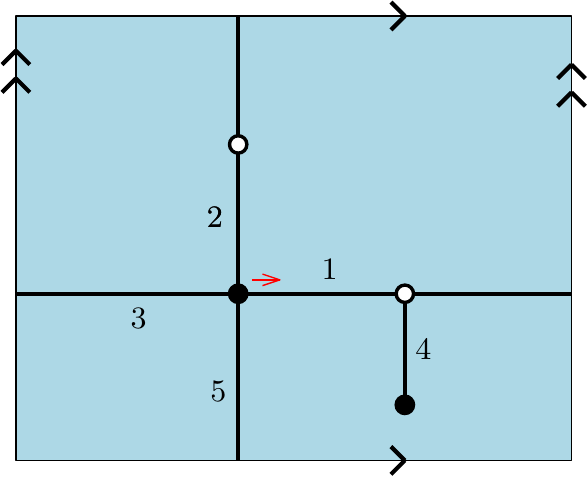}
\caption{}
\label{fig:bipartite map}
\end{subfigure}
\caption{(A) A rooted map corresponding to the permutations $\nu=(1\,8\,6\,13)(2\,9\,4\,14\,11)(3\,7\,10\,12)(5)$, $\epsilon=(1\,7)(2\,14)(3\,5)(4\,13)(6\,12)(8\,11)(9\,10)$, $\phi=\epsilon \nu = (1\,11\,14\,8\,12\,5\,3)(2\,10\,6\,4)(7\,9\,13)$. (B) A rooted bipartite map corresponding to the permutations $\nu_{\bullet}=(1\,5\,3\,2)(4)$, $\nu_{\circ}=(1\,3\,4)(2\,5)$ with $\phi=\nu_{\circ}\nu_{\bullet}=(1\,2\,3\,4\,5)$. Both maps are embedded in the torus.}
\label{fig:maps}
\end{figure}

\begin{itemize}
    \item We decorate a rooted map with $n$ edges by assigning labels from $1,2,\ldots,2n$ to the two ends of each edge with the restriction that end of the root edge corresponding to the vertex has label $1$. To such a decorated rooted map we associate three permutations $\nu, \epsilon, \phi$ in $\Sym_{2n}$ as follows: (i) the clockwise list of ends of each vertex give the cycles of the {\em vertex permutation} $\nu$, (ii) the pairs of labels on each edge give the transpositions of the {\em edge permutation} $\epsilon$ (a fixed-point-free involution), and (iii) the counterclockwise list of the second label on each edge when traversing each face gives the cycles of the {\em face permutation} $\phi$.  By construction we have $\epsilon \nu = \phi$. See Figure~\ref{fig:map}. 
    \item If the rooted map with $n$ edges is bipartite then the vertices can be properly coloured by two colours, say black and white, and we assume the root is coloured black. We assign the labels $1,2,\ldots,n$ to each of the edges with the restriction that the root edge has label $1$. To such decorated rooted bipartite map we associate three permutations $\nu_{\bullet}, \nu_{\circ}, \phi$ in $\Sym_n$ as follows: (i) the clockwise list of edge labels incident to each black vertex give the cycles of the {\em $\bullet$-vertex permutation} $\nu_{\bullet}$, (ii) the clockwise list of edge labels incident to each white vertex give the cycles of the {\em $\circ$-vertex permutation} $\nu_{\circ}$, (iii) the counterclockwise list of the label on each edge $\bullet-\circ$ when traversing each face gives the cycles of the {\em face permutation} $\phi$. Again, by construction we have $\nu_{\circ}\nu_{\bullet} = \phi$. See Figure~\ref{fig:bipartite map}. 
\end{itemize}
The genus of the map can be read from the number of cycles of the associated triple of permutations by Euler's formula. The transitivity of $\langle \nu, \epsilon \rangle$ and  $\langle \nu_{\bullet}, \nu_{\circ} \rangle$ on the labels $[2n]$ and $[n]$, respectively, follows from the connectivity of the respective maps.

In the early 90s, David Jackson with his former Ph.D. student Terry Visentin \cite{jacksonvisentin,jacksonvisentin2} studied rooted maps using this connection to factorizations of permutations and the machinery of irreducible characters and the group algebra of the symmetric group. With this approach they found in \cite[Cor. 5.2]{jacksonvisentin} the following remarkable identity between rooted maps and quadrangulations. Let $M(u,x,y,z)$ and $Q(u,x,y,z)$ be the generating series for rooted maps and rooted quadrangulations where the indeterminates $u,x,y,z$ mark the genus, number of faces, vertices, and edges, respectively. Then
\begin{equation} \label{eq: rel maps and quadrangulations}
2 Q(u^2,x,y,z) = M(4u^2,y+u,y,xz^2) +  M(4u^2,y-u,y,xz^2). 
\end{equation}
In \cite{JPV}, Jackson, Visentin and the physicist Perry used this identity to  prove a conjectured connection between two models in physics at the level of {\em Feynman diagrams}. Jackson and Visentin asked for a bijective proof of this result. This problem is known as the {\em Quadrangulation conjecture}. 

\begin{conjecture}[Quadrangulation conjecture] \label{prob: quad conj}
There is an explicit bijection that preserves genus and number of edges between the set of rooted quadrangulations and rooted maps with certain decorations that proves \eqref{eq: rel maps and quadrangulations}.
\end{conjecture}

The genus zero case of such a bijection is the {\em medial bijection} of Tutte \cite{Tutte63}. Other special cases towards finding such a bijection are due to Jackson and his former Ph.D. student D.R.L. Brown \cite{BJ}.

\subsection{Transitive factorizations into transpositions and Hurwitz numbers} \label{sec:Hurwitz}

A third constraint that can be placed on factors is \emph{global}.
Referring to \eqref{eq:notrans}, we can additionally require
\begin{equation}\label{eq:transdef}
	c)  \textnormal{ the factors } \tau_1, \dots, \tau_k \textnormal{ generate a
	transitive subgroup of } \Sym_n. \hspace{3cm}
\end{equation}
A \emph{transitive factorization} of $\rho$ satisfies all three criteria of \eqref{eq:notrans} and \eqref{eq:transdef}.

Let $\alpha  = (\alpha_1, \ldots, \alpha_{l(\alpha)}) \vdash n$ be a partition of $n$ with $l(\alpha)$ parts and let $\rho \in \conj_\alpha$.    An interesting case of the transitive factorization problem is when the factors are all transpositions.  We refer to this as \textbf{Problem 1}.  It can be shown that a transitive factorization of $\rho$ requires at least $n+l(\alpha)-2$ transpositions;  parity considerations then dictate that the number of transitive factorizations of $\rho$ is nonzero only if the number of factors is $n+l(\alpha)-2 + 2g$ for some $g$.  We call $g$ the \emph{genus} of a factorization, with the case $g=0$ being referred to as \emph{minimal} transitive factorizations.  By symmetry considerations, we see that the number of such factorizations of $\rho$ is only dependent on its conjugacy class,  and we let $b^g_\alpha$ be the number of transitive factorizations of $\rho$ with $n+l(\alpha)-2 + 2g$ factors.  Remarkably, the numbers $\bee^g_\alpha$ are related to the number of topologically inequivalent $n$-sheeted branched covers of genus $g$ of the sphere, where the branching type is given by $\alpha$ at one specified point, and the other  branching points, totalling $n+l(\alpha) - 2 + 2g$, are given by transpositions.  The numbers of such branched covers are known as the \emph{Hurwitz numbers} of genus $g$ and type $\alpha$ and are denoted by  $H^g_\alpha$.  The Hurwitz numbers and transitive factorizations into transpositions satisfy the relationship $H^g_\alpha = \tfrac{|\conj_\alpha|}{n!} \bee^g_\alpha$.  Goulden and Jackson themselves have written a survey \cite{goulsurvey} of transitive factorizations of permutations, and the extensive connections between these objects and geometry.

\subsubsection{Single Hurwitz numbers}
A remarkably simple formula exists for  $\bee^0_\alpha$ attributed to Hurwitz:  
\begin{equation}\label{eq:bee}
    \bee^0_\alpha = {n^{l(\alpha)-3}} (n+l(\alpha)-2)! \prod_{i=1}^{l(\alpha)} \frac{\alpha_i^{\alpha_i}}{(\alpha_i-1) !}.
\end{equation}
Hurwitz \cite{hurwitz} gave a sketch of a proof by induction, reconstructued in detail by Strehl \cite{strehl}. Goulden and Jackson proved this formula through a \emph{join-cut} analysis.  The analysis involves focusing on the effect of the last transposition $\tau_k = (i\; j)$ in a product $\tau_1 \cdots \tau_k$ on the product of the previous $k-1$ factors.  Considering the product $\hat{\rho} = \tau_1 \cdots \tau_{k-1}$:
\begin{itemize}
    \item if $i$ and $j$ are in the same cycle of $\hat{\rho}$, then they will be in separate cycles in $\tau_1 \cdots \tau_k$ ($\tau_k$ \emph{cuts} a cycle of $\hat{\rho}$); or
    \item if $i$ and $j$ are in different cycles of $\hat{\rho}$, then they will be in the same cycle in $\tau_1 \cdots \tau_k$ ($\tau_k$ \emph{joins} cycles of $\hat{\rho}$).
\end{itemize}

To prove \eqref{eq:bee}, a generating series for $\bee^0_\alpha$ is constructed, and the join-cut analysis is then used to show this generating series satisfies a partial differential eqution whose solution is unique.  An ingenious use of Lagrange Inversion is then used to show that the generating series for the numbers on the right hand side of \eqref{eq:bee} also satisfy the differential equation;  whence the equality in \eqref{eq:bee} follows.  The details are found in \cite{goul:2}.  Goulden and Jackson later discovered that the formula in \eqref{eq:bee} was known to Hurwitz \cite{hurwitz}. Later, Bousquet-M\'elou and Schaeffer \cite{B-M--S} gave a bijective construction using maps to show a generalization of \eqref{eq:bee} allowing factors that are not transpositions. See also \cite{DPSI,DPSII} for another bijective approach to Hurwitz numbers.

Goulden and Jackson \cite{GJ99} and Vakil \cite{V01} also proved the following formula for $b^1_{\alpha}$, conjectured by Goulden--Jackson--Vainshtein \cite{GJVain}:
\[
b^1_{\alpha} \,=\, \frac{1}{24} (n+m)! \left(\prod_{i=1}^m
  \frac{\alpha_i^{\alpha_i}}{(\alpha_i-1)!}\right) \left( n^m -
  n^{m-1}  - \sum_{i=2}^m (i-2)! \cdot e_i(\alpha)\cdot  n^{m-i}\right),
\]
where $e_i(\alpha)$ denotes the $i$th elementary symmetric function evaluated at $\alpha=(\alpha_1,\ldots,\alpha_m)$.

Lastly, there has been work to generalize Hurwitz numbers to other groups. Bini, Goulden, and Jackson in [BGJ08] consider extensions to the hyperoctahedral group, and recently in \cite{PR} and \cite{DLM1,DLM2}, the authors look at generalizations in complex reflection groups with different notions of transitivity.

\subsubsection{Double Hurwitz numbers} 

Another important transitive factorization problem of a $\rho \in \conj_\alpha$ is to require all but one factors to be transpositions, and the remaining factor is from a specified conjugacy class $\conj_\beta$.  We refer to this as \textbf{Problem 2}.  In this case it can be shown that the number of transpositions required is $l(\alpha) + l(\beta) + 2g -2$, where again $g$ is referred to as the genus of such factorizations.  We label the number of such factorizations of $\rho$ by $\bee^g_{\alpha, \beta}$.  The numbers $\bee^g_{\alpha, \beta}$ are related to the \emph{double Hurwitz numbers} $H^g_{\alpha, \beta}$, which were introduced by Okounkov \cite{OkKP} in his study of Gromov-Witten theory. These numbers count the number of topologically inequivalent $n$-sheeted branched covers of genus $g$ with two specified branch points of type $\alpha$ and $\beta$ and the other branching points are given by transpositions.  Several of these aspects of transitive factorizations of permutations into transpositions are discussed in Goulden and Jackson's own survey \cite{goulsurvey}, including their work with Vakil \cite{GJV1,GJV2} on the polynomiality and quasi-polynomiality of single and double Hurwitz numbers, respectively, using {\em Ehrhart theory} of polytopes (see also \cite{CJM,CMbook}).  We discuss here an aspect not mentioned there connected to work of one of the survey authors (AR).

\subsubsection{Star factorizations}

Let $i \leq n$ be positive integers.  Let $S_i$ be the set of transpositions in
$\Sym_n$ containing the symbol $i$.  Thus $|S_i| = n-1$.  Let $\alpha =
(\alpha_1, \ldots, \alpha_m)$ be a partition of $n$ and $\conj_\alpha$ be
the associated conjugacy class in $\Sym_n$.   Consider the problem of writing a permutation
$\rho \in \conj_\alpha$ as the transitive product of transpositions in $S_i$.  We refer to this factorization problem as \textbf{Problem 3}.
The set $S_i$ is intimately connected with the \emph{Jucys--Murphy} elements,
famously used by Okounkov and Vershik \cite{okver} in their study of the
representation theory of the symmetric group.  Like Problem 1 for the single
Hurwitz numbers, the number of
transpositions needed for a transitive factorization of $\rho$ into factors
coming from $S_i$ is $n+l(\alpha) + 2g -2$ for a nonnegative number $g$
called the genus of the factorization.     \emph{A priori} the number of factorizations of $\rho$
of fixed length should depend on the length of the cycle containing the symbol
$i$ in $\rho$ as $i$ is distinguished in these factorizations;  however, David Jackson's former Ph.D.
student John Irving and 
Amarpreet Rattan
discovered in \cite{rattanirving:1} that number of such factorizations for
$g=0$ (the minimal case) is
\begin{equation}\label{eq:starsym}
 \frac{(n+m-2)!}{n!} \cdot \alpha_1\cdots \alpha_m.
\end{equation}
Note the symmetry in \eqref{eq:starsym};  that is, the cycle containing the
symbol $i$ in $\rho$ plays no special role in the formula.  Subsequently, Goulden and Jackson showed that
this symmetry persists for $g > 0$;  that is, the number of such
factorizations in genus $g > 0$ is only dependent on the conjugacy class
$\alpha$ of the target permutation $\rho$ and not the specific permutation
itself.  Accordingly, let $a^g_{\alpha}$ be the number of transitive
factorizations of a fixed permutation in the conjugacy class $C_\alpha$ of
genus $g$.  Thus $a^0_\alpha$ is given in \eqref{eq:starsym}, while Goulden and
Jackson give a formula for $a^g_\alpha$ for $g > 0$.

When the number of factorizations of a permutation $\rho$ in a transitive
factorization problem is only dependent on the conjugacy class of $\rho$, we call
the problem \emph{central}.  Thus Problems 1 and 2, which pertain to the single
and double Hurwitz numbers, respectively, are obviously central through
symmetry considerations, while Problem 3 is central through the results in
\cite{rattanirving:1, goul:8}.  Intriguingly, the methods in both
\cite{rattanirving:1, goul:8} show the centrality of Problem 3 as corollaries of
the full enumerative formulas for $a^g_{\alpha}$, and they do not expose
combinatorially why Problem 3 is central.  Combinatorial
constructions showing the centrality of Problem 3 in the $g=0$ case have been
subsequently found \cite{tenner, tenner2}, but the general case is still open.

Also in \cite{goul:8}, Goulden and Jackson show a compact relationship between
the numbers $a^g_\alpha$ and the double Hurwitz numbers $b^g_{\alpha, \beta}$.  For a partition $\alpha$ of $n$ with $m=l(\alpha)$, they show
\begin{equation}\label{eq:transnotrans}
    a^g_{\alpha} = \frac{1}{n! (2n-1)^{n+m-3+2g}} b^g_{\alpha \cup 1^{n-1}, (2n-1)}.
\end{equation}
Here, when $\alpha$ is a partition of $n$, the partition $\alpha \cup 1^{n-1}$ is the partition of $2n-1$ obtained from $\alpha$ by inserting $n-1$ parts of size 1.
Thus \eqref{eq:transnotrans} directly connects Problem 3 to the obviously central Problem 2.
Again, the proof by Goulden and Jackson for \eqref{eq:transnotrans} is to compare the final enumerative
results of the numbers $a^g_\alpha$ and the double Hurwitz numbers rather than
giving a combinatorial construction connecting the two.  Goulden and
Jackson ask for a direct combinatorial explanation of \eqref{eq:transnotrans},
still open for all $g \geq 0$. 

\subsubsection{Monotone Hurwitz numbers}

We also discuss a variant of Hurwitz numbers that was studied in a series of papers \cite{GG-PN1,GG-PN2,GG-PN3,GG-PN4,GG-PN5} by Ian Goulden, his former Ph.D. student Mathieu  Guay-Paquet, and his former postdoc Jonathan Novak. 

A variation of the transitive factorization problem is to require the transpositions $\tau_i=(a_i\,b_i)$ with $a_i<b_i$ in the factorization to have the restriction that $b_1\leq b_2 \leq \cdots \leq b_k$. Such factorizations are called {\em monotone} and are related to expansions of complete symmetric functions in the Jucys--Murphy elements \cite{MN}. These monotone factorizations have nice enumerative formulas. For instance, Daniele Gewurz and Francesca Merola \cite{GerMer} showed that out of the $b^0_{(n)} = n^{n-2}$ factorizations of the long cycle in $\Sym_n$ into $n$ transpositions, there are Catalan many $\frac{1}{n+1}\binom{2n}{n}$ monotone factorizations. Thus, in general, referring to \eqref{eq:notrans} and \eqref{eq:transdef}, for monotone factorizations we additionally require
\begin{equation}
\label{eq:monotone} d) \textnormal{ writing each factor } \tau_i=(a_i,b_i) \textnormal{ with } a_i<b_i, \textnormal{ we have } b_1\leq \cdots \leq b_k. 
\end{equation}
We refer to this factorization problem as \textbf{Problem 4} and denote by $\overrightarrow{b}_{\alpha}^g$ the number of transitive monotone factorizations of $\rho$ with $n+l(\alpha)-2+2g$ factors. In \cite{GG-PN2}, the authors used a modified join-cut approach to prove the following formula for genus zero monotone factorizations
\begin{equation} \label{eq:genus zero monotone Hurwitz numbers}
\overrightarrow{b}_{\alpha}^0 \,=\, n!\cdot (2n+1)^{\overline{l(\alpha)-3}}\cdot \prod_{j=1}^{l(\alpha)} \binom{2\alpha_j}{\alpha_j},
\end{equation}
where $(m)^{\overline{k}}=m(m+1)\cdots (m+k-1)$, denotes a rising factorial. Note the striking similarity between this formula and Hurwitz's formula \eqref{eq:bee} for $b^0_{\alpha}$. See \cite[Thm. 2.1]{BerIrv} for an elegant reciprocity relation between regular factorizations and monotone factorizations; and see \cite{BCD} for recent connections between monotone factorizations and integrable hierarchies (see Section~\ref{sec:KP}).

\subsection{KP hierarchy} \label{sec:KP}

In \cite{GJ_KP}, Goulden and Jackson brought to algebraic combinatorics a connection between  {\em integrable systems} and factorizations/maps. The {\em Kadomtsev--Petriashvili (KP) hierarchy} is an infinite list of partial differential equations for a function $F=F(p_1,p_2,\ldots)$ that generalizes the famous {\em Korteweg--De Vries (KdV) equations} that model waves in shallow water. The first two equations of this hierarchy are
\begin{align} \label{eq:KP1}
F_{2,2} - F_{3,1} +\frac{1}{12} F_{1^4} + \frac12 (F_{1^2})^2 &= 0\\
\intertext{and}
\label{eq:KP2}
F_{3,2} - F_{4,1} + \frac16 F_{2,1^3} + F_{1,1}F_{2,1} &=0,
\end{align}
where $F_{r^{a_r},\ldots,1^{a_1}} := \dfrac{\partial^{a_r+\cdots + a_1}}{\partial p_r^{a_r}\cdots \partial p_1^{a_1}} F$.

Goulden and Jackson showed in \cite{GJ_KP} that the generating functions of transitive factorizations are solutions of the KP hierarchy. Given partitions $\alpha$ and $\beta$ of $n\geq 1$ and $a_1,a_2,\ldots \geq 0$, let  $b_{\alpha,\beta}^{(a_1,a_2,\ldots)}$ be the number of ordered tuples $(\sigma,\gamma,\pi_1,\pi_2, \ldots)$ of permutations on $\mathfrak{S}_n$ such that 
\begin{equation}\label{eq:notrans double}
\begin{aligned}
	a)\;\;\;\;\;\;&  \sigma \gamma \pi_1 \pi_2 \cdots = \iota \text{ where $\iota$ is the identity}, \\
	b)\;\;\;\;\;\;&  \sigma \in  \conj_{\alpha}, \gamma \in \conj_{\beta}, \text{ and } n - l(\pi_i)=a_i  \text{ for }   i=1,\ldots,k\\
	c)\;\;\;\;\;\;& \textnormal{the factors } \sigma, \gamma, \pi_1, \dots, \pi_k \textnormal{ generate a	transitive subgroup of } \Sym_n. 
\end{aligned}
\hspace{1cm}
\end{equation}
 Note that this number generalizes the number of factorizations  $\gamma(a_1,\ldots,a_p;n)$ from \eqref{eq:jackpfactors} and the numbers $b^g_{\alpha}$, $b^g_{\alpha,\beta}$ related to the single (double) Hurwitz numbers mentioned in Section~\ref{sec:Hurwitz}. For instance
\begin{align*}
b^g_{\alpha} &=  b^{(1^r,0,\ldots)}_{\alpha,1^n} \quad \text{for}\quad  r=l(\alpha)+n+2g-2.
\end{align*}

There is a classical characterization of solutions to the KP hierarchy that are logarithms of linear combinations of Schur functions:
\(\sum_{\lambda} b_{\lambda} s_{\lambda}(p_1,p_2,\ldots)\),
where $b_{\lambda} \in \mathbb{Q}[u_1,u_2,\ldots]$  whenever the $b_{\lambda}$ satisfy the famous {\em Pl\"ucker relations} from the study of the {\em Grassmannian}. This characterization combined with a result of Orlov and Shcherbin \cite{OS}, imply that certain linear combinations of Schur functions, involving {\em contents} of Young diagrams, satisfy the KP hierarchy. Goulden and Jackson noticed that since the generating functions of the numbers  $b_{\alpha,\beta}^{(a_1,a_2,\ldots)}$ of factorizations can be written in terms of such Schur functions involving contents, then these generating functions  are also solutions to the KP hierarchy. 

\begin{theorem}[Goulden--Jackson \cite{GJ_KP}] \label{thm:GJ KP}
The series 
\[
\sum_{\substack{|\alpha|=|\beta|=n \geq 1 \\ a_1,a_2,\ldots \geq 0}} \frac{1}{n!} b_{\alpha,\beta}^{(a_1,a_2,\ldots)} p_{\alpha} q_{\beta} u_1^{a_1}u_2^{a_2}\cdots, 
\]
is a solution of the KP hierarchy (in the variables $p_1,p_2,\ldots$).
\end{theorem}

As applications of this result, they show that the KP hierarchy is satisfied by the generating series for (i)  double Hurwitz numbers (recovering a result of Okounkov \cite{OkKP}), (ii) the number of rooted maps on orientable surfaces, and (iii) the number of triangulations in an orientable surface. 

Soon after Carrell  (a former Ph.D. student of Goulden) and Chapuy studied in \cite{CarC} the {\em first} equation \eqref{eq:KP1} of the KP hierarchy applied to the number of factorizations/maps to derive a new quadratic recurrence for the number of orientable rooted maps of a given genus  that dramatically improved the efficiency of counting such maps (see \cite[\S 5.5]{Scha_survey}). 

This connection between integrality and enumeration of factorizations/maps continues to be a very active and rich avenue of research  (see \cite{DYZ,BCD}).

\subsection{The \texorpdfstring{$b$}--conjecture and the matchings-Jack conjecture}
\label{sec: b conjecture}

Let $s_{\lambda}({\mathbf x})$ and $p_{\lambda}({\mathbf x})$ denote the {\em Schur} and {\em power sum} symmetric functions indexed by $\lambda$ in the variables ${\mathbf x} = (x_1,x_2,\ldots)$. The connection coefficients $c^{\gamma}_{\alpha,\beta}$ of the symmetric group defined in \eqref{eq:concoeff}
appear in the following elegant identity of {\em symmetric functions}. 
\begin{equation} \label{partition function S_n}
\sum_{n\geq 1} t^n \sum_{\alpha,\beta,\gamma \vdash n}    c^{\gamma}_{\alpha,\beta} \frac{|\conj_\gamma|}{n!} p_{\alpha}({\bf x})p_{\beta}({\bf y}) p_{\gamma}({\bf z}) ,=\, \sum_{\lambda} \hooks(\lambda) s_{\lambda}({\bf x}) s_{\lambda}({\bf y})s_{\lambda}({\bf z}) t^{|\lambda|}, 
\end{equation}
where  the sum on the right-hand-side is over all  partitions $\lambda$, and $\hooks(\lambda) = \prod_{(i,j) \in \lambda} (\lambda_i+\lambda'_j-i-j+1)$ is the product of the hook lengths of $\lambda$. 

In \cite{GJ-Jack}, Goulden and Jackson started the study of a generalization of \eqref{partition function S_n} in terms of {\em Jack symmetric functions} $J^{(\alpha)}_{\lambda}({\bf x})$. These symmetric functions specialize to Schur functions for $\alpha=1$, to {\em Zonal polynomials} for $\alpha=2$ and are related to {\em Macdonald polynomials} (see \cite{StJ,Mc}). Goulden and Jackson defined the rational functions $c^{\gamma}_{\alpha,\beta}(b)$ as follows.
\begin{equation} \label{Jack gs}
 \sum_{n\geq 1} t^n \sum_{\alpha,\beta,\gamma \vdash n}    \frac{c^{\gamma}_{\alpha,\beta}(\alpha-1)}{\alpha^{l(\gamma)}} \frac{|\conj_\gamma|}{n!}p_{\alpha}({\bf x})p_{\beta}({\bf y}) p_{\gamma}({\bf z}) = \sum_{\lambda} \frac{1}{\hooks_{\alpha}(\lambda)\hooks'_{\alpha}(\lambda)} J^{\alpha}_{\lambda}({\bf x})J^{\alpha}_{\lambda}({\bf y})J^{\alpha}_{\lambda}({\bf z}) t^{|\lambda|},
\end{equation}
where $\hooks_{\alpha}(\lambda)$ and $\hooks'_{\alpha}(\lambda)$ are certain simple $\alpha$-deformations of $\hooks(\lambda)$ (see \cite[Thm. 5.8]{StJ}). Note that at $b=\alpha-1=0$ we recover the connection coefficients of the symmetric group $c^{\gamma}_{\alpha,\beta}(0)=c^{\gamma}_{\alpha,\beta}$. Goulden and Jackson in \cite{GJ-Jack} have the following polynomiality and integrality conjecture for $c^{\gamma}_{\alpha,\beta}(b)$.

\begin{conjecture}[Matchings-Jack conjecture]
For all partitions $\lambda,\mu,\nu \vdash n$ the connection coefficient $c^{\gamma}_{\alpha,\beta}(b)$ is a polynomial in $\mathbb{N}[b]$. Moreover, there is a statistic $\stat_{\gamma}$ on certain matchings $\delta$ of $2n$ points such that
\[
c^{\gamma}_{\alpha,\beta}(b) = \sum_{\delta} b^{\stat_{\gamma}(\delta)}.
\]
\end{conjecture}

Goulden and Jackson also defined the rational functions $h^{\lambda}_{\mu,\nu}(b)$ as follows.
\begin{equation}
\sum_{n \geq 1} t^n \sum_{\alpha,\beta,\gamma \vdash n} h^{\alpha}_{\beta,\gamma}(b) p_{\alpha}({\bf x})  p_{\beta}({\bf y})  p_{\gamma}({\bf z}) \,=\, (b+1)t\frac{\partial}{\partial t} \log \Phi({\bf x},{\bf y},{\bf z};t,b+1),   
\end{equation}
where $\Phi({\bf x},{\bf y},{\bf z};t,\alpha)$ is the generating function in \eqref{Jack gs}. Again, at $b=\alpha-1=0$ and $\beta=\alpha-1=1$ the quantity $h^{\lambda}_{\mu,\nu}(b)$ counts certain maps in orientable and  locally orientable surfaces. Goulden and Jackson in \cite{GJ-Jack} have the following similar polynomiality and integrality conjecture for $h^{\gamma}_{\alpha,\beta}(b)$.

\begin{conjecture}[$b$-conjecture]
For all partitions $\lambda,\mu,\nu \vdash n$ the connection
coefficient $h^{\lambda}_{\mu,\nu}(b)$ is a polynomial in
$\mathbb{N}[b]$. Moreover, there is a statistic $\stat$ on locally orientable maps
such that
\[
h^{\lambda}_{\mu,\nu}(b) = \sum_{M} b^{\stat(M)},
\]
where the sum is over all rooted, bipartite locally orientable maps $M$ with face
distribution $\lambda$, black vertex distribution $\mu$, and white
vertex distribution $\nu$. Moreover $\stat(M)=0$ if and only if $M$ is orientable.
\end{conjecture}

These two conjectures are according to \cite{CD} ``among the most remarkable open problems in algebraic combinatorics'' since for $b\neq \{0,1\}$ there are very few tools to approach them. The conjecture is related to representation theory of {\em Gelfand pairs} \cite{HSS}, random matrices, and algebraic geometry. Special cases of the conjecture have been proved in \cite{Lacroix} by  M.~A. La Croix, a former Ph.D. student of Goulden and Jackson, and in \cite{BJ,KV,KPVass}. Moreover, there has been important recent progress on these conjectures: polynomiality over $\mathbb{Q}[b]$ of both conjectures was proved by F\'eray--Do{\l}{\k{e}}ga in \cite{DF1,DF2}, the case ${\bf z} = (z,z,\ldots)$ (i.e. keeping track of two sets of variables ${\bf x}$ and ${\bf y}$) was settled by  Chapuy--Do{\l}{\k{e}}ga \cite{CD}, and polynomiality over $\mathbb{Z}[b]$ for the first conjecture was settled by Ben Dali \cite{dali2022integrality} using the {\em Farahat--Higman} algebra \cite{FH}.

\section{Towards an algebraic foundation for quantum field theory}
\label{quantum}

Over a series of three papers \cite{1404.0747, MR3659109, Jackson_2021}, following a program outlined in an earlier preprint \cite{0810.42931}, David Jackson, together with collaborators Achim Kempf and Alejandro Morales worked to give an algebraic foundation to quantum field theory.  They extended this into practical algorithms for integration which fell out of the field theory work in \cite{Kempf_2015}.

Quantum field theory is the quantum theory of particles interacting.  The term quantum field theory is used generically to refer to all such theories, while a specific quantum field theory is determined by the particular particles and interactions being considered.  The standard model in fundamental particle physics is a specific quantum field theory and is highly successful, giving some of the most precise calculations of experimentally measured values anywhere in science and describing everything we currently understand in high energy physics.  However, quantum field theory suffers from many problems in its mathematical foundation. In particular, quantum field theories are often described in terms of a path integral, but the path integral needs to be an integral over fields, and a satisfactory analytic definition remains a question of active research.

Jackson, Kempf and Morales are setting the foundations for a different approach using formal power series to give an algebraic foundation to quantum field theory.

Other researchers, including one of the survey authors (KY) \cite{Ybook},  work in this direction.
The work of Jackson, Kempf and Morales stands apart for its focus on transforms and its intentional approach to the question of foundations for quantum field theory.

They summarize the picture of quantum field theory as
\begin{equation}\label{eq picture}
e^{iS[\Phi]} \overset{Fourier}{\leftrightarrow} Z[J] \overset{log/exp}{\leftrightarrow} iW[J] \overset{Legendre}{\leftrightarrow} i\Gamma[\phi]
\end{equation}
where $S[\Phi]$ is the action, $Z[J]$ is the path integral, which perturbatively can be interpreted as a sum over all graphs, while $iW[J]$ is the sum over connected graphs and $i\Gamma[\phi]$ is the sum over 1 particle irreducible graphs, that is, connected bridgeless graphs.

The justification usually given for using these transforms is based on analytic conditions such as convexity, but we know from the theory of algebraic enumeration that the exponential map takes us from connected objects to potentially disconnected objects -- this is the \emph{combinatorial exponential map} -- and it makes sense on any formal power series with zero constant term.  In a similar way, Jackson, Kempf and Morales define combinatorial Legendre \cite{MR3659109} and Fourier \cite{Jackson_2021} transforms making all steps of \eqref{eq picture} combinatorial under appropriate hypotheses.

The proof of the combinatorial Legendre transform is essentially Euler's formula for trees $1=V(g)-E(g)$ for a graph $g$, interpreted appropriately, and so the combinatorial Legendre transform is the ``tree-of'' operator on combinatorial objects. This transform is also related to the Lagrange Inversion theorem. The combinatorial Fourier transform in a similar way is the ``graph-of'' operator.

As a predecessor to the general work on the combinatorial Fourier transform, Jackson, Kempf and Morales took a formal power series look at the Dirac delta in \cite{Kempf_2015}.

In particular, for a sufficiently nice function $g$, the Dirac delta can be written as
$\delta(x) = \frac{1}{\sqrt{2\pi}}\frac{1}{g(-i\partial_x)} \widetilde{g}(x)$,
where $\widetilde{g}$ is the Fourier transform of $g$.
{}From this and similar equations
they gave new approaches to integrating by differentiating \cite{Kempf_2015}, which were subsequently further developed by Kempf and collaborators \cite{JTK} and used in computer algebra systems for practical computations.

Together this collection of papers gives an algebraic and combinatorial foundation to key transforms of quantum field theory and moves us towards a full formal foundation for quantum field theory.

\section{Acknowledgements}

We thank  Theo Douvropoulos, Ian Goulden, and Joel Brewster Lewis for comments and suggestions. We also thank Vic Reiner, Dennis Stanton, and Sheila Sundaram with help with historical background. This work was supported by the Canadian Tri-Council Research Support Fund. The author A.M.F. was supported by an NSERC Discovery Grant. K.Y. was supported by an NSERC Discovery Grant and the Canada Research Chairs program.  A.H.M. was partially supported by the NSF grants DMS-1855536 and DMS-22030407. 

\bibliographystyle{amsplain-ac}
\bibliography{combined_refs}

\end{document}